\documentclass{article}
\usepackage{amsmath}
\usepackage{amsfonts}
\usepackage{amssymb}

\begin{document}

%%%%%%%%%%%%%%%%%%%%%%Definitions%%%%%%%%%%%%%%%%%%%%%%%%%%%%%%%%%%%%%%%%%%%

\renewcommand{\theequation}{\thesection.\arabic{equation}}

\newcommand{\hsp}{{\hspace*{\parindent}}}
\newtheorem{theorem}{Theorem}[section]
\newtheorem{problem}{Problem}[section]
\newtheorem{definition}{Definition}[section]
\newtheorem{lemma}{Lemma}[section]
\newtheorem{proposition}{Proposition}[section]

\newtheorem{corollary}{Corollary}
\newtheorem{example}{Example}
\newtheorem{conjecture}{Conjecture}
\newtheorem{algorithm}{Algorithm}
\newtheorem{exercise}{Exercise}
\newtheorem{remarkk}{Remark}

\newcommand{\req}[1]{(\ref{#1})}

\newcommand{\lip}{\langle}
\newcommand{\rip}{\rangle}
\newcommand{\uu}{\underline}
\newcommand{\oo}{\overline}
\newcommand{\La}{\Lambda}
\newcommand{\la}{\lambda}
\newcommand{\eps}{\varepsilon}
\newcommand{\Om}{\Omega}
\newcommand{\om}{w}

\newcommand{\dint}{\displaystyle\int}
\newcommand{\dsum}{\displaystyle\sum}
\newcommand{\dfr}{\displaystyle\frac}
\newcommand{\bige}{\mbox{\Large\it e}}
\newcommand{\integers}{{\mathbb Z}}
\newcommand{\ZZ}{{\Bbb Z}}
\newcommand{\rationals}{{\mathbb Q}}
\newcommand{\reals}{{\mathbb R}}
\newcommand{\realsd}{\reals^d}
\newcommand{\realsn}{\reals^n}
\newcommand{\NN}{{\mathbb N}}
\newcommand{\DD}{{\mathbb D}}
\newcommand{\XX}{{\mathfrak X}}
\newcommand{\calA}{\mathcal A}
\newcommand{\calC}{\mathcal C}
\newcommand{\calL}{\mathcal L}
\newcommand{\calF}{\mathcal F}
\newcommand{\calP}{\mathcal P}
\newcommand{\calX}{\mathcal X}

\newcommand{\dfn}{\stackrel{\triangle}{=}}
\def\complex{\mathop{\raise .45ex\hbox{${\bf\scriptstyle{|}}$}
     \kern -0.40em {\rm \textstyle{C}}}\nolimits}
\def\hilbert{\mathop{\raise .21ex\hbox{$\bigcirc$}}\kern -1.005em {\rm\textstyle{H}}} %Hilbert space
\newcommand{\RAISE}{{\:\raisebox{.6ex}{$\scriptstyle{>}$}\raisebox{-.3ex}
           {$\scriptstyle{\!\!\!\!\!<}\:$}}} % >< one above each other

\newcommand{\pp}{{\partial}}
\newcommand{\al}{{\alpha}}
\newcommand{\bfcdot}{{\mbox{\boldmath$\cdot$}}}
\newcommand{\bolda}{\mathbf{a}}
\newcommand{\boldb}{\mathbf{b}}
\newcommand{\boldq}{\mathbf{q}}

\newcommand{\dett}{{\textstyle{\det_2}}}
\newcommand{\sign}{{\mbox{\rm sign}\,}}
\newcommand{\TE}{{\rm TE}}
\newcommand{\TA}{{\rm TA}}
\newcommand{\E}{{\rm E\,}}
\newcommand{\won}{{\mbox{\bf 1}}}
\newcommand{\Lebn}{{\rm Leb}_n}
\newcommand{\Prob}{{\rm Prob\,}}
\newcommand{\sinc}{{\rm sinc\,}}
\newcommand{\ctg}{{\rm ctg\,}}
\newcommand{\loc}{{\rm loc}}
\newcommand{\trace}{{\,\rm trace\,}}
\newcommand{\Tracenew}{{\,\,\rm trace\,}}
\newcommand{\Dom}{{\rm Dom}\,}
\newcommand{\ifff}{\mbox{\ if and only if\ }}
\newcommand{\proof}{\noindent {\bf Proof:\ }}
\newcommand{\remark}{\noindent {\bf Remark:\ }}
\newcommand{\remarks}{\noindent {\bf Remarks:\ }}
\newcommand{\note}{\noindent {\bf Note:\ }}

\newcommand{\one}{\frac{1}{n}\:}
\newcommand{\half}{\frac{1}{2}\:}
\newcommand{\limn}{\lim_{n \rightarrow \infty}}
\newcommand{\W}{\raisebox{-.25cm}{$_W$}}
%qed
\def\squarebox#1{\hbox to #1{\hfill\vbox to #1{\vfill}}}
\newcommand{\qed}{\hspace*{\fill}
           \vbox{\hrule\hbox{\vrule\squarebox{.667em}\vrule}\hrule}\bigskip}
%%%%%%%%%%%%%%%%%%%%%%%%%%%%%%%%%%%%%%%%%%%%%%%%%%%%%%%%%%%%%%%%%%%%%%%%%%%%
\newcommand{\ls}[1]
   {\dimen0=\fontdimen6\the\font \lineskip=#1\dimen0
\advance\lineskip.5\fontdimen5\the\font \advance\lineskip-\dimen0
\lineskiplimit=.9\lineskip \baselineskip=\lineskip
\advance\baselineskip\dimen0 \normallineskip\lineskip
\normallineskiplimit\lineskiplimit \normalbaselineskip\baselineskip
\ignorespaces }

%%%%%%%%%%%%%%%%%%%%%%%%%%%%%%%%%%%%%%%%%%%%%%%%%%%%%%%%%%%%%%%%%%%%%%%%%%%%

\title{Rotations and Tangent Processes on Wiener Space}
\author{M. Zakai}

\date{}

\maketitle

%\thispagestyle{empty}
%\setcounter{page}{0}

%\vspace{2cm}
\ls{1.5}
\begin{abstract}
\noindent
The paper considers (a) Representations of measure preserving transformations
(``rotations'') on Wiener space, and (b) The stochastic calculus of variations
induced by parameterized rotations $\{T_\theta w, 0 \le \theta \le \eps\}$:
``Directional derivatives'' $(dF(T_\theta w)/d \theta)_{\theta=0}$,
``vector fields'' or ``tangent processes'' 
$(dT_\theta w /d\theta)_{\theta=0}$ and flows of rotations.
\end{abstract}

\section{Introduction}
\hsp
Let $(W, H, \mu)$ be an abstract Wiener space (AWS): $W=\{w\}$ is a separable 
Banach
space, $H$ (the Cameron-Martin space) is a separable Hilbert space  
densely and continuously embedded in $W$,
$W^*\hookrightarrow H^* = H \hookrightarrow W$ and for every $e$ in $W^*$,
$_{w^*}(e,w)_w$ is $N(0,|e|^2_H)$.
By the Cameron-Martin theorem, for any $h\in H$, the measure induced by $w+h$
is equivalent to the measure $\mu$, therefore if $F(w)$ is a r.v.\ on the
Wiener space, so is $F(w+h)$.  This fact
enabled the development of the stochastic calculus of variations, i.e. the
Malliavin calculus which very roughly is based on the directional derivative
of $F$ in the $h$ direction:  
$(dF(w+\eps h)/d\eps)_{\eps=0}$.  Now, let $T$ be a measure preserving
transformation on $W$ (in short, a `Rotation'), i.e.
$_{w^*}(e, Tw)_w$ is also
$N(0,|e|^2_H)$.  Then if $F(w)$ is a r.v.\ so is $F(Tw)$ and if 
$T_\theta w$ $0\le \theta \le \eps$ is a smooth collection of rotations one 
can
consider objects like $(d F(T_\theta w)/d\theta)_{\theta=0}$.
The purpose of this paper is to survey previous work and to present new
results on  the following:
\begin{itemize}
\item[(i)]
Measure preserving transformations (rotations)
\item[(ii)]
The Malliavin calculus of rotations.
\end{itemize}

The study 
of stochastic analysis
over Riemannian manifolds showed that the Cameron-Martin space is not
sufficient to represent the tangent space and 
as discussed in
(\cite{D92,FM,D95,CM96,M,CC}),  
more general vector fields are
needed. 
The setup of these papers was based on
the model of
\begin{equation}
\label{1.1}
(Tw)_i = \int_0^\bfcdot \sum_{j=1}^d \sigma_{ij} (s,w) dw_j(s), 
\qquad i,j\le d
\end{equation}
for ``rotational vector fields'',
where $w$ is the d-dimensional $(d \ge 2)$ Wiener process with $\sigma$ 
being non-anticipative 
and
skew symmetric.  This paper considers the abstract Wiener space setup 
presented in \cite{HUZ},
it is restricted to flat space.
A particular class of anticipative tangent processes was recently 
considered in \cite{CM00}.  A class of rotations on Wiener space 
introduced in \cite{K2} are different from the rotations considered here.

Section~2 summarizes results, mainly from the Malliavin Calculus, which are
needed later.
Rotations $T$ are considered in section~3, where $Tw=\sum \eta_i (w) e_i, $ $
\eta_i(w)$ are i.i.d.\ $N(0,1)$ random variables, and $e_i$,
$i=1,2,\cdots$ is a complete orthonormal base in $H$. 
Rotations are introduced in Section~3.
Theorem 3.1 (\cite{UZ95,UZ99}) presents rotations by showing that sequences of i.i.d.\ $(N(0,1)$ random variables
can be constructed as the divergence of $R(w) e_i$ where
$R(w)$ belongs to a certain class of operators.
An outline of the proof is included.  This is followed by new results,
Proposition~3.1 and Theorem~3.2
showing that under some smoothness assumptions every sequence of i.i.d.\
$N(0,1)$ random variables on the Wiener space can be represented by the
construction of Theorem~3.1.
Section~4 deals with directional derivations of the type $(dF(T_\theta
w)/d\theta)_{\theta=0}$.  A ``tangent operator'' is introduced and  its
relation to the directional derivative is indicated.  Section~5 deals with
``tangent process'' which are Banach valued random variables, play the role
of tangent vectors and induce the directional derivatives.  The
first part of section~6 gives a  positive answer (due to Tsirelson and
Glasner) to a problem raised in \cite{HUZ} whether the group  of invertible
rotations on the Wiener space are connected. The second part deals with flows
of rotations, i.e. the flow induced by
$$
\frac{d T_tw}{dt} = m(T_t w, t)
$$
where $m(w,t)$ are the tangent processes introduced in section~5.  
The case where $m(w,t)$ is of the type of equation of \req{1.1} was considered
in \cite{CC}.
A more
detailed proof of the result of \cite{HUZ} is given.  The appendix deals
with the following problem:  In view of the results of section~3 and other
results, the question arose whether the condition that $\nabla u(w)$  
 is quasinilpotent implies the existence of a filtration
such that $u$ is adapted.  A counter example, following \cite{R} is presented
in the appendix.

\paragraph{Acknowledgements:} We
wish to express our sincere thanks to B.\ Tsirelson and E. Glasner for their
contribution in Section~6, to A.S. \"Ust\"unel for useful discussions and to
E. Mayer-Wolf for useful remarks.
\section{Preliminaries}
\setcounter{equation}{0}

\paragraph{Notation:}
For each $e\in H$ and induced by an element of $W^*$, $\delta (e)$ denotes the
$N(0, |e|_H^2)$ r.v.\
$_{W^*}(e,w)_W$.  For all $e\in H$, $\delta (e)$ denotes the $L^2$ limit of 
$_{W^*}(e_n, w)_W$ as $e_n \to e$ in $H$.  

We will not distinguish between embeddings and inclusions.  For example, 
$e\in W^*$ will also be considered as an element of $H$ or $W$; the distinction
being clear from context.

For $F(w) = f(\delta e_1, \cdots,
\delta e_n)$ and $f$ smooth, $\nabla F$ is defined as
$\nabla F = \sum$ $ f_i'(\delta e_1, \cdots , \delta e_n) \cdot e_i$.  Let
$\calX$ be a separable Hilbert space and $u$ an $\calX$ valued functional
$\DD_{p,k} (\calX)$, $p>1$, $k \in  \NN$ will denote the Sobolev space of
$\calX$ valued functionals in $L^P(\mu, \calX)$ whose $k$-th order derivative
$\nabla^k u$ is in $L^P(\mu, \calX \otimes H^{\otimes k})$.
$\DD_{p,k}(R)$ will be denoted $\DD_{p,k}$.
$\DD(\calX) = \cap_{p>1} \cap_{k\in N} \DD_{p,k} (\calX)$.  Recall that
$$
\nabla: \: \DD_{p,k} (\calX) \to \DD_{p,k-1} (\calX \otimes H)
$$
and for $\delta$, the adjoint of $\nabla$ under the Wiener measure
$$
\delta:\: \DD_{p,k} (\calX \otimes H) \to \DD_{p,k-1} (\calX)
$$
are continuous linear operators for any $p>1, k \in \NN$.
The operator $\delta$ is the divergence  or the Skorohod integral and:
\begin{itemize}
\item[(a)]
If $u \in \DD_{2,1} (H)$, then
$$
E[(\delta u)^2] = E[|u|_H^2] + E [\trace (\nabla u)^2]\,.
$$
\item[(b)]
If $F\in \DD_{2,1}$, $u \in \DD_{2,1}(H)$ and if $Fu \in \DD_{2,1} (H)$, then
\begin{equation}
\label{2.1}
\delta (Fu) = F\delta u - (\nabla F, u)_H\,.
\end{equation}
\end{itemize}

\paragraph{A. Exact and divergence free $H$ valued r.v's.:}\mbox{}\\ 
Let $u\in \DD_{2,1} (H)$ then (a)~ $u$ is said to be
``exact'' if $u\in \nabla F(w)$ for some $\DD_{2,1}$ functional
$F(w)$.
(b)~$u$ is said to be divergence
free if $\delta u=0$.

Set
\begin{align*}
U^e =
 U^{\text{exact}} & = \{ u\in \DD_{2,1} (H) : u=\nabla F\} \\
U^{\text{d.f.}} =
U^{\text{diverg. free}} & = \{ u \in \DD_{2,1} (H) : \delta u=0\}
\end{align*}

If $u\in U^e, v \in U^{\text{d.f.}}$ then
$E(u,v)_H = E(\nabla f, v)_H= E(f \delta v) = 0$.  Hence 
$U^e$ and $U^{\text{d.f.}}$ are orthogonal subspace of $\DD_{2,1}(H)$
and $\DD_{2,1} (H) = U^e \oplus  U^{\text{d.f.}}$. 

Let $\calL F$ be the Ornstein-Uhlenbeck operator:
$\calL F = \delta \nabla F$, assume that $EF=0$ then $\calL^{-1}$
is a bounded operator and $\calL \calL^{-1} F=F$.  Hence, for
any
$F\in \DD_{2,1}, EF=0 F(w)= \delta (\nabla \calL^{-1} F)$ and $F$ possesses the
representation $F(w) = \delta u, u = \nabla \calL^{-1} F$, where
$u\in \DD_{2,1} (H)$.  Note that this representation is different from the
Ito-type representation of Wiener functionals.

Returning to $U^e$ and $U^{\mathrm{d.f.}}$, for 
$F(w) = \delta u, u \in \DD_{2,1} (H)$,
$\delta (u-\nabla \calL^{-1} \delta u) = 0$ then
$\nabla \calL^{-1} \delta u$ and $u-\nabla \calL^{-1} \delta u$ are the
projections of $u$ on $U^e$ and $U^{\mathrm{d.f.}}$ respectively.

We prepare, for later reference, the following lemma.
\begin{lemma}
Let $u\in \DD_{2,1} (H)$, let $\{e_i, i=1, 2, \cdots \}$ be a CONB on $H$
further assume that $u=\sum \delta (v_i) e_i,\: v_i \in \DD_{2,1} (H)$.
If $ (v_i, e_j)_H + (v_j, e_i)_H = 0, \; i,j=1, 2, \cdots$,
then $\delta u=0$.
In particular the above result holds for
 $v_i = A(w) e_i$, 
where $A+A^T=0$. 
\end{lemma}

\proof
For smooth $F$, integrating by 
parts we have
\begin{align*}
%\label{newno}
E(F\cdot \delta u)   
& = E \sum_i \nabla_{e_i} F \delta v_i\notag \\
& = E \sum_i \nabla_{e_i} F \cdot \delta \left(\sum_j (v_i, e_j)
e_j\right)\notag \\
& = E \sum_{i,j} \nabla_{e_i, e_j}^2 F \cdot (v_i, e_j)\,.
\end{align*}
and $\delta u=0$ follows
since $\nabla_{e_i, e_j}^2 F$ is symmetric in $i$ and $j$ and $F$ is arbitrary.

\paragraph{B.\ Constructing a filtration on the AWS}\mbox{}\\
Let $(W,H,\mu)$ be an A.W.S.,  we introduce
a time structure i.e.\ a filtration and causality on it as follows:
Let the projections
$\{\pi_\theta, \; 0 \le \theta \le 1\}$ be a continuous strictly increasing
resolution of the identity on  $H$.  Set
$$
\calF_\theta = \sigma \{ \delta \pi_\theta h, \quad h \in H\}
\,.
$$
Propositions \ref{prop-2.1}--\ref{prop-2.3} are from \cite{UZ97}.

\begin{proposition}
\label{prop-2.1}
$F(w) \in \DD_{2,1}$ is $\calF_\theta$ measurable iff
$\nabla F = \pi_\theta \nabla F$ a.s.
(intuitively: if $F(w) = f(\delta h_1, \cdots \delta h_n)$ with
$h_i = \pi_\theta h_i$)
\end{proposition}
\begin{definition}
An  $H$-valued r.v. $u$ will be said to be $\calF_\bfcdot$-adapted
if for every $\theta \in [0,1]$ and $h\in H$, $(u, \pi_\theta h)$ is
$\calF_\theta$ measurable.
\end{definition}
\begin{proposition}
\label{prop-2.2}
$u\in \DD_{2,1} (H)$ is $\calF_\bfcdot$ adapted iff
$$
\pi_\theta \nabla u \pi_\theta = \pi_\theta \nabla u
$$
for all $\theta \in [0,1]$
\end{proposition}
\begin{definition}
Let $G(w)$ be a measurable r.v.\ taking values in the class of
bounded transformations on $H$.  Then $G$ is said to be \textit{weakly
adapted\/} if for all $h\in H$, $Gh$ is adapted, $G$ will be said to be
\textit{adapted\/} (or causal) if $Gu$ is adapted for all adapted $u$.
\end{definition}
\begin{proposition}
\label{prop-2.3}
If $G(w)$ satisfies $Gu \in \DD_{2,0} (H)$ whenever $u\in \DD_{2,0}(H)$ and
$G$ is weakly adapted then $G$ is adapted iff
$\pi_\theta G\pi_\theta = \pi_\theta G$.
\end{proposition}
Another version of the last result is:
\begin{proposition}
Under the assumptions of the previous proposition, a weakly adapted 
$G$ is adapted iff for all $h\in H$
\begin{equation}
\label{star}
\pi_\theta h =0 \Rightarrow \pi_\theta Gh=0
\,.
\end{equation}
\end{proposition}

\proof
Let $u$ be of the form
\begin{equation}
\label{starstar}
u=\sum_1^n \varphi_i \Bigl( \pi_{\theta_{i+1}} - \pi_{\theta_i} \Bigr) h_i
\end{equation}
where $\theta_{i+1} > \theta_i$ and the $\varphi_i$ are $\DD_2$ r.v.'s.
Then
\begin{equation}
\label{starstarstar}
A u = \sum_1^n \varphi_i G \Bigl( \pi_{\theta_{i+1}} - \pi_{\theta_i}\Bigr)
h_i
\end{equation}
$h_i \in H$.  Now, assume that $u$ is adapted hence the $\varphi_i$ are 
$\calF_{\theta_i} $ measurable.  Since $G$ is weakly measurable then \req{star}
implies that $\varphi_i G \pi_{\theta_{i+1}} (I - \pi_{\theta_i}) h_i$ is also
adapted 
hence $Gu$ is adapted and $G$ is adapted
since $u$ of the form \req{starstar} are dense in $\DD_{2} (H)$.

Conversely, again $u$ is assumed to be adapted, and since $Gu$ is adapted,
$\varphi_i G (\pi_{\theta_{i+1}} - \pi_{\theta_i}) h_i$ is adapted.  Hence,
since $\varphi_i$ are $\calF_{\theta_i}$ measurable, we must have
$$
\pi_\theta G \Bigl( \pi_{\theta_{i+1}} - \pi_{\theta_i}\Bigr) h_i = 0
$$
and \req{star} follows.
\qed

Given a $\DD_{2,0}$ functional $F(w)$ on $(W, H, \mu)$ and a filtration
$\pi_\theta$ (continuous strictly increasing) then
there exists a unique adapted $u \in \DD_2 (H)$ such that 
$u\in \Dom \delta$, and
$F(w) = \delta u$,
and $E(\delta u)^2 = E(|u|_H^2$.
(\cite{UZ00},
in the classical setup this follows directly 
from the multiple Wiener integral).
Hence, given  $F(w) = \delta u^e$
then $u^e$ can be
lifted uniquely to $\tilde u$ such that $\delta \tilde u = \delta u^e$ and
$\tilde u$ is adapted to a given filtration.  

\paragraph{C.\ Quasinilpotent operators}\mbox{}\\
An H-S operator on $H$ is said to be quasinilpotent (q.n.p.) if any
\textit{one} of the following is satisfied (cf. \cite{ZZ} or \cite{UZ99}):
\begin{itemize}
\item[(a)]
$\trace A^n = 0 \qquad \forall n \ge 2$.
\item[(b)]
$|A^n |^{\one} \longrightarrow 0, \qquad |\cdot |$ is the operator norm.
\item[(c)]
The spectrum of $A$ is $\{0\}$ only.
\item[(d)]
$(1-\alpha A)^{-1} = \sum_{n=0}^\infty \alpha^n A^n \qquad \forall \alpha$.
\item[(e)]
$\det_2 (I+\alpha A) = 1\quad \forall \alpha$, where
$ \det_2 (I+A) = \Pi (1-\la_i) e^{-\la_i}$.
\end{itemize}
\begin{proposition}[\cite{UZ97}]
If $u$ is adapted, $u \in \DD_{2,1} (H)$ then $\nabla u$ is q.n.p.\
\end{proposition}
\paragraph{Outline of proof:} 
Let $\theta_{i+1} > \theta_i$, set
$$
u_i = q_i \Bigl( \pi_{\theta_{i+1}} - \pi_{\theta_i} \Bigr)
h_i = q_i \tilde{h}_i\,.
$$
Then,
\begin{align*}
\trace \nabla u_i \nabla u_j & =
\trace (\nabla q_i \otimes \tilde h_i) (\nabla q_j \otimes \tilde{h}_j) \\
& = (\nabla q_i, \tilde h_j)_H (\nabla q_j, \tilde h_i)_H\,.
\end{align*}
For $i=j$,
$(\nabla q_i, \tilde h_j)_H =0$, if $i>j$ then
$(\nabla q_i, \tilde h_i)_H = 0$.  Similarly for
$i<j$.

The following question arises regarding the converse of the last result:
Given $u$ such that $\nabla u$ is q.n.p., does this imply the  existence of 
a filtration such
that $u$ is adapted to it?  The answer to this question, 
as shown in appendix A,  is negative. 

\paragraph{D.\ The  Ito-Nisio theorem}
\begin{theorem}
Let $(X_i)$ be a symmetric sequence of random variables (i.e.\linebreak 
$(\pm X_1, \pm X_2, \cdots \pm X_n)$ has the same law as
$(X_1, \cdots X_n)$ for any $n$)
with values in a separable Banach space $B$.  Denote by $\mu_n$ the
distribution of the partial sum $S_n=\sum_{i=1}^n X_i$.  The following are
equivalent:
\begin{itemize}
\item[(i)]
the sequence $(S_n)$ converges almost surely in the Banach norm;
\item[(ii)]
$(S_n)$ converges in probability;
\item[(iii)]
$(\mu_n)$ converges weakly;
\item[(iv)]
there exists a $B$-valued r.v.\ $\gamma$ such that $(S_n,f)
\underset{P}{\longrightarrow} (\gamma, f)$ for all $f $ in $B'$;
\item[(v)]
there exists a probability measure $\mu$ in $\calP(B)$ sucht that
$\mu_n \circ f^{-1} \to \mu \circ f^{-1}$ weakly for every $f$ in $B'$.
\end{itemize}
\end{theorem}
Cf.\ \cite{LT} for (i), (ii), (iii) and (v); (iv) follows from (i) and implies (v).

Let $(W, H, \mu)$ be AWS, then for any CONS, $\{e_i\}$
and any i.i.d. $N(0,1)$ random variables $\eta_i$ then
$$y_n = \sum_1^n \eta _i e_i
$$
converges a.s.\ in the Banach norm
(in particular $|w-\sum_1^n \delta e_i e_i|_W
\overset{\text{a.s.}}{\longrightarrow} 0$).
Hence denoting $y=\lim y_n $, then 
$y=Tw$ is a measure preserving transformation and 
$\eta_i =\,_W(Tw, e_i)_{W^*}$.

\section{Rotations on Wiener space}
\setcounter{equation}{0}

\hsp
Let $\eta_i = \text{i.i.d.} \quad N(0,1)$ r.v.'s and $\{e_n\} $ a 
CONB induced by $W^*$
then by the Ito-Nisio theorem
\begin{equation}
\label{3.1}
Tw=\sum \eta_i e_i
\end{equation}
is a measure preserving transformation, we will refer to it as a rotation.

%\bigskip
%\noindent
%{\bf Theorem 3.1 \cite{UZ95,UZ99} -- The Rotation Theorem.}
%{\em 
\begin{theorem}[\cite{UZ95,UZ99}]
Let $w\mapsto R(w)$ be a strongly measurable random variable on $W$ with values in the
space of bounded linear operators on $H$.  Assume that $R$ is almost surely 
an
isometry on $H$, i.e. $|R(w) h|_H$ $ = |h|_H$ a.s.\ for all
$h\in H)$.  Further assume that for some $p>1$ and for all 
$h\in H$,
$Rh \in \DD_{p,2}(H)$, and $\nabla Rh \in \DD_{p,1} (H\otimes H)$
is a quasi-nilpotent operator on $H$.   If moreover, either

\noindent
(a)~$(I+i\nabla Rh)^{-1} \cdot Rh$ is in $L^q(\mu,H), q > 1$ for any
$h\in H$ (here $q$ may depend on $h\in H$) or,\\
(b)~$Rh \in \DD(H)$ for a dense set in $H$.

Then
$$
%\label{2.1new}
E\Bigl[\exp i \delta (Rh)\Bigr] = \exp - \half |h|^2_H\,.
$$

\noindent
Namely, if $(e_n, n\in\NN)$ is a complete, orthonormal basis in $H$
then
$(\delta(R e_n), n\in\NN)$ are independent 
$N(0,1)$-random variables and consequently $\sum_i \delta(R e_i) e_i$ defines
a measure preserving transformation of $W$.
%}
\end{theorem}

The map $ R$ satisfying the conditions for this theorem with $p=2$ and under (a)
with $q=2$ will be said to
{\em satisfy the rotation conditions}.

\paragraph{Outline of proof:} Let $u: B\to H$ be ``an H-C$^1$ map'' and
$T=w+u$, assume that $T$ is a.s.\ invertible then \cite{UZ99}.
$$
E\Bigl(F(Tw) \cdot | \La (w) | \Bigr) = EF(w)
$$
where
$$
\La(w) = \dett (I_H + \nabla u) \exp \left( - \delta u - \half |u|_H^2\right)
\,.
$$
In particular, $F(w)=1$,
$u=\nabla Rh$, then since $\det_2(\qquad) = 1$ and
$|Rh|_H^2 = |h|_H^2$, hence
$$
E\exp \left( -\delta (Rh) - \half 
|h|_H^2\right) =1
$$
or
$$
E\exp - \delta (Rh) = \exp - \half |h|_H^2
$$
consequently $\delta(Re_i), i=1,2, \cdots $ are i.i.d., $N(0,1)$ and
$\sum \delta (R e_i) e_i $ is a rotation by the Ito-Nisio theorem.

The conditions on $R$ in the notation theorem are obviously not necessary
since if for all $i$, 
$u_i \in U^{\text{d.f.}}$ and defining
$\rho:H\to H$ by
$$
\rho e_i = u_i, \qquad i=1,2,\cdots
$$
then $\delta \rho h=0$ and if $R$ induces a rotation so does $R+\rho$.  We
have, however, the following two results, which yield
a converse to the rotation theorem.
\begin{proposition}
Let $R(w)$ be an a.s.\ bounded operator on $H$.   Assume that 
$R(w)$ is weakly adapted with respect to a filtration induced by a continuous
increasing $\pi_\bfcdot$.  Further assume
for all
$h\in H$, $R(w) h$ is in the domain of $\delta$ and the probability law of
$\delta(Rh)$ is $N(0, |h|_H^2)$ then:
\begin{enumerate}
\item
If $h_1, h_2 \in H$ and $(h_1, h_2)_H=0$ then $\delta(Rh_1)$ and
$\delta(Rh_2)$
are independent.
\item
$R(w)$ is a.s.\ an isometry on $H$.
\item
$\sum_i \delta(R e_i) e_i$ is a rotation and if $e_i, i=1, 2, \cdots$ and
$h_i, i=1, 2, \cdots$ are CONB's on $H$ then, a.s.
$$
\sum_i \delta (R h_i) h_i = \sum_i \delta (R e_i) e_i\,.
$$
\end{enumerate}
\end{proposition}

\proof 
\begin{align*}
1.\quad E\Bigl(\exp\, i \alpha\, \delta (R h_1) \exp\, i \beta\, \delta (R h_2)\Bigr)
& = E\exp\, i \delta \Bigl(R(\alpha h + \beta h_2)\Bigr) \\
& = \exp - \frac{\alpha^2}{2} |h_1|_H^2 - \frac{\beta^2}{2} |h_2|_H^2 \\
& = E \Bigl( \exp\, i \alpha \, \delta Rh_1\Bigr) E \exp \, i \beta \, \delta (R h_2)
\end{align*}
2.\quad By part 1,
$y_\theta = \delta (R \pi_\theta h)$ is a Gaussian process of independent
increments.  Hence it is Gaussian martingale and  its quadratic variation
satisfies
$$
\lip y, y \rip_\theta = E y_\theta^2 
\,.
$$
and by our assumption $Ey_\theta^2 = |\pi_\theta h|_H^2$.
But
$$
\lip y, y \rip_\theta = (R \pi_\theta h, R \pi_\theta h)_H
$$
and $R^T R=I$ follows.

\noindent
3.\ Follows from the Ito-Nisio theorem.
\qed

\setcounter{theorem}{1}
\begin{theorem}
Let $Tw =\sum \eta_i e_i$ be a rotation.  Let $\calF_\bfcdot$ be a filtration
induced by a continuous increasing resolution of the identity.  Then there
exists a unique $\calF_\bfcdot$ weakly adapted $R(w): H\to H$ which is an isometry and
$$
Tw=\sum_i \delta (Re_i) e_i\,.
$$
If, moreover, $\eta_i \in \DD_{2,2}$ then $\nabla Rh$ is q.n.p.
\end{theorem}

\proof
By our assumptions, every $\eta_i$ can be uniquely represented as
$\eta_i=\delta u_i$ where the $u_i$ are adapted, in the domain of $\delta$,
and $u_i \in \DD_{2} (H)$.  Define $R$ by
$$
R(w) e_i = u_i
$$
then $R(w)$ is weakly adapted, and satisfies the assumptions of 
the previous result.  Hence $R$ is an isometry and
$Tw=\sum\delta (Re_i) e_i$.  If moreover the $\eta_i \in \DD_{2,2}$ then
$\nabla Rh \in \DD_{2,1} (H)$ and is q.n.p. since it is adapted.
\qed

\remark
For a given rotation $Tw=\sum \delta (R e_i)e_i$, 
$R(w)$ is ``highly non unique'';
instead of representing $\eta_i$ as the divergence of adapted processes, we
can define $\eta_i =\delta v_i$ with $v_i \in U^e$ to yield a  unique $R^e(w)$
such that $R^eh$ is exact for all $h\in H$ and $Tw = \sum\delta (R^ee_i) e_i$.
In other words, given any $R$ satisfying the assumptions of the theorem we can
construct an $R^e$ such that $\delta(R^eh) = \delta(Rh)$ and
$R^eh\in U^e$ for all $h\in H$.  Thus $R^e$ will not necessarily be an
isometry. Also, we can ``lift'' $R^e$ to another $\tilde R$ which
is weakly adapted with respect to another filtration.

\paragraph{Two examples of rotations:}\mbox{}\\
(a)~Let $R$ be a rotation.  Assume that $R(w)$ and 
$R^T(w)$ are adapted (not just weakly adapted) and $R^T(w) = R^{-1}(w)$.  Then
$$
R(w) = \int_0^1 \exp 2 \pi i \theta\, \pi (d\theta, w)
$$
and \cite[section~109]{RN} $\pi(\theta,w)$ can be approximated by
polynomials in $R$ and $R^T$.  Consequently for 
$\{\theta \in [0,1]), \pi(\theta, w)\}$
is
adapted.  Therefore
$$
A(w) = \log R(w) = \int_0^1 2 \pi i \theta \, \pi (d\theta, w)
$$
is adapted $A+A^T =0$, for $\alpha \in \reals$,
$$
R^\alpha (w) = \int_0^1 \exp 2 \pi i \alpha \theta  \,\pi(d\theta, w)
$$
is also adapted
and under some additional conditions 
$R^\alpha (w)$ induces a rotation for all real $\alpha$.

\noindent
(b)~Let $\pi(\theta)$ be a nonrandom resolution of the identity and
$\calF_\bfcdot$ the filtration induced by it.  Let $C$ denote the class of
unitary transformations on $H$ which commute with $\pi(\theta)$ for every
$\theta$.  Let
$$
R(w) = \int_0^1 B(\theta, w) \pi (d\theta)
$$
and assume that $B(\theta,w)\in C$ and 
$B(\theta,w) \pi(\theta) h, \theta \in b [0,1]$, $h\in H$, is
$\calF_\theta$ measurable.  Then $R(w)$ is weakly adapted and unitary so that
under some technical condition  $R$ is a rotation.  Moreover by our
assumptions on $B$, $\pi(\theta) R = R\pi(\theta), \theta \in [0,1]$ hence
$R(w)$ is weakly adapted.  Similarly, let $D$ denote the class of skew symmetric
operators $(K+K^T = 0)$ which commute with $\pi(\theta)$.  Set
$$
A(w) = \int_0^1 b (\theta, w) \pi (d\theta)
$$
where $b\in D$ and $b(\theta, w) \pi(\theta) \cdot h$ is $\calF_\theta$ 
measurable.  Then
$A+A^T=0$, $A$ is weakly adapted and
$\int_0^1 \exp b(\theta,w) \pi(d\theta)$ is a rotation.
Note that \req{1.1} is a special case of this example.

\section{Tangent operators}
\setcounter{equation}{0}

\hsp
Let $T_tw=\sum_i \delta (R_t e_i) e_i$ where $R_t, t \in [0,
\delta]$ is unitary, satisfies the rotation condition,
and let $(e_i, i=1, 2, \cdots)$ be a CONB on $H$ induced by $W^*$.
 Assume
that $(R_{t+\eps} - R_t)/\eps$ converges a.s.\ in the operator norm to
$B_t (w)$ as $\eps \to0$, then 
$(R_{t+\eps}^T - R_t^T)/\eps$ converges to $B_t^T$ and 
$0=\frac{d(R_t^T R_t)}{dt} = \frac{dR_t R_t^T}{dt}$.
Hence
\begin{equation}
\label{4.1}
B_t^T R_t + R_t^T B_t  = 0
\mbox{\ and \ }
B_t R_t^T + R_t B_t^T  =  0\,.
\end{equation}
Setting $R_0 = I$ and $A=(d R_t/dt) $ at $t=0$ yields
\begin{equation}
\label{4.2}
A+A^T = 0\,.
\end{equation}
Let $f(x_1, \cdots, x_n)$ be a smooth function $\reals_n$.  Set $F(w) =
f(\delta e_1, \cdots, \delta e_n)$.  Since
$_W(T_tw, e_i)_{W^*} = \delta (R_t e_i)$ and $\delta e_i \circ T_t = 
_W\!\!(T_t w, e_i)_{W^*}$,
$\delta e_i \circ T_\eps w = _W\!\!(T_tw, e_i)_{W^*}$.

Hence, with $A=(d R_t/dt)_{t=0}$:
\begin{align*}
\left.  \frac{dF(T_t w)}{dt}\right|_{t=0}
& = \sum_i \left(f_i' \Bigl(\delta (R_te_1), \cdots, \delta (R_te_n)\Bigr)
\delta
\left(\frac{dR_t}{dt} e_i\right)\right)_{t=0} \\
& = \sum_i f_i' (\delta e_1, \cdots, \delta e_n) \delta (Ae_i)\\
& = \sum_i \delta \Bigl(f_i' (\cdot) Ae_i\Bigr) - \sum_{ij} f_{ij}'' (\cdot) (e_j, Ae_i)
\end{align*}
The second term vanished since $f_{ij}''$ is symmetric and $A$ is skew
symmetric.  Hence
\begin{equation}
\label{new4.3}
\left.\frac{dF(T_t w)}{dt}\right|_{t=0} = \delta (A \nabla F)
\end{equation}

Motivated by \req{new4.3} 
we define
\begin{definition}[\cite{HUZ}]
Let $w\to Q(w)$ be a weakly measurable mapping taking values in the space of
bounded operators on $H$.  Assume that for all $F\in \DD_{2,1}$,
$Q\nabla F \in \DD_{2,1} (H)$.  For every $Q$ satisfying these conditions and
$u\in \DD_{2,1} (H)$ we define
$$
\calL_{Q,u} F = \delta (Q \nabla F) + \nabla_u F
$$
and denote it as the Tangent Operator induced by $(Q, u)$.  Also
$\calL_{Q,0} F =: \calL_Q F$.  The following summarizes some properties of
the tangent operator (cf.\ \cite{HUZ} for proofs).
\end{definition}
\begin{itemize}
\item[1)]
$\calL_{Q,u}$ is closeable in $H$ (i.e. if $F_n\to 0$ in $H$ a.s.\ and $\calL
F_n$ exist, then $\calL F_n\to 0$).
\item[2)]
The adjoint of $\calL_{Q,u}$ satisfies $\calL_{Q,u}^* F = \calL_{Q} F +
\delta (Fu)$.

\item[3)]
If $Q=A$ where
$A^T+ A = 0$, then
$$
\calL_A F_1 F_2 = F_1 \calL_A F_2 + F_2 \calL_A F_1
$$
namely, $\calL_A$ is a derivation (i.e.\ behaves as a first order operator).
\item[4)]
For $F\in \DD$,  $e^F \nabla F \in \DD_{2,1} (H) $
and $A+A^T=0$
$$
\calL_A e^F = e^F \calL_A F + (\nabla F, A \nabla F) = e^F \calL_A F
$$
and for all $f\in C_b^1$
$$
\calL_A f(\delta h) = f' (\delta h) \delta (A h)
\,.
$$
\item[5)]
Cf.\ \cite{HUZ} for results for $\calL_A (\delta u)$ and $[\calL_A,
\calL_B]$.
\end{itemize}

Let $R$ be unitary and satisfy the rotation condition.  Let $R_{t,k}(w)$ denote
$R(w+ t\cdot k)$, $t\in [0,1]$, $k\in H$.  Then $R_{t,k}$
is also a.s.\ unitary and $\nabla R_{t,k} h$ is also a.s. quasinilpotent.
Assume that $R_{t,k}$ satisfies condition $a$ or $b$ of the rotation theorem
then $R_{t,k}$ also induces a rotation, let $T_{t,k}$ denote this rotation.
Setting
\begin{align*}
X_k^R F(w) & = \left. \frac{dF(T_{t,k} w)}{dt} \right|_{t=0} \\
& = \calL_{\dot{R}_k} F
\end{align*}
where
$$
\dot{R}_k = \left. \frac{dR_{t,k}}{dt} \right|_{t=0}
$$
%By equation \req{4.1} $R^{-1}(w) X_k^R$ and 
%$X_k^R R^{-1}$ are skew symmetric, and 
As shown in \cite{UZ95} and \cite{UZ99}
\begin{equation}
\label{4.3}
\nabla (F \circ T) = R(\nabla F \circ T) + X^R F
\end{equation}
i.e.\
$$
\nabla_h (F \circ T) = \Bigl( R(\nabla F \circ T), h\Bigr)_H + X_h^R F
$$
and when $W\to H$ is a cylindrical map
\begin{align}
\label{4.4}
(\delta u) \circ T &= \delta \Bigl(R (u\circ T)\Bigr) + 
\sum_i \Bigl(X^R(u,e_i), Re_i\Bigr)\notag\\
 &= \delta \Bigl(R (u\circ T)\Bigr) + \trace (R^{-1} X^R u)\,.
\end{align}
If $F(w)$ and $g(x), x\in \reals$ is smooth then
\begin{align*}
0 & = \frac{d}{dt} E g\Bigl( F (T_{t,k} w) \Bigr) \\
  & = E \left\{ g' \Bigl( F(T_{t,k} (w) )\Bigr) \cdot X_{t,k}^R F
\right\}
\,.
\end{align*}

\section{Tangent processes}
\setcounter{equation}{0}
\hsp
Let $R_t$ satisfy the rotation condition, set
$ \calL_{R_t} w = T_tw = \sum_i \delta (R_t e_i ) e_i  $.
In order to represent the ``directional derivative'' $\calL_A F$ as the action
of a ``tangent vector'' on $F$, the ``vector field'' $(d T_tw/dt)_{t=0}$ is
needed.  Formally, for $A=\left(\frac{d R_t}{dt}\right)_{t=0}$,
$$
\left.\frac{dT_tw}{dt}\right|_{t=0} = \sum_i \delta (A e_i) e_i
$$
which motivates the following definition:

\begin{definition}
\label{tan-def}
Let $Q(w)$ be a weakly measurable $H$ operator valued
transformation on $H$.  Assume that $Q(w)h \in
\DD_{2,1}(H)$ for all $h\in H$. 
Let $e_i, i=1,2,\cdots$ be a CONB induced by elements of $W^*$, if
$\sum_i \delta (Q e_i)
 e_i$ converges weakly  in the Banach space
as $n \to \infty$; 
namely, if there exists a 
$W$-valued random variable $Y$ such that
$\sum_i \delta(Q e_i) (e_i, \tilde\alpha)_H$ 
converges in probability 
($\tilde\alpha$ is the image of $\alpha$ in $H$ under
the canonical injection from $W^*$ to $H$) to $_W\lip Y,\alpha\rip_{W^*}$
for all $\alpha \in W^*$.  The limit $Y$
will be denoted $Y= \calL_Q w$ and will be
called the tangent process induced by
$Q$ and $\{e_i\})$.
\end{definition}

\remarks
(a)~The definition given here is somewhat different from that in \cite{HUZ}
as $Y$ may depend on $\{e_i\}$.

(b)~It is often necessary to consider the case where the series 
$\sum\delta (Q e_i) e_i$ satisfy a stronger convergence condition; several
cases assuring a.s.\ convergence are
\begin{itemize}
\item[a.] The case where the $\delta (Qe_i)$ satisfy the conditions
of the extended Ito-Nisio theorem.  

\item[b.] The case where $Q$ is a bounded non random operator
(\cite[Theorem~1.14]{K}).
\item[c.] 
Let $\varphi \in \DD_{2,1}$ and $Q^* \nabla \varphi \in \DD_{2,0} (H)$, then
$\calL_{\varphi Q} w$ exists in the sense
of a.s.\ convergence in $W$ if and only if $\varphi \calL_Q w$ exists in the
corresponding sense and then
$$
\calL_{\varphi Q} w = \varphi \calL_Q w - Q^T \nabla \varphi\,.
$$
\end{itemize}
The proof follows directly from
$$
\sum_1^N \delta(\varphi Qe_i) e_i = \sum_1^N \varphi\delta(Q_i e_i) e_i-
\sum_1^N (Q^T \nabla \varphi, e_i) \cdot e_i\,.
$$
The relation between the tangent process $\calL_Q w$ and the tangent
operator is reflected in the following lemma.
\begin{lemma}
Assume that $Q$ satisfies the requirements of Definition~\ref{tan-def}
and\linebreak 
$\sum \delta (Qe_i)e_i$ converges a.s.  Further assume that
$u\in \DD_{2,1} (H), Qu \in \DD_{2,1} (H)$ $u$ is the image in $H$ of 
$\tilde{\tilde u}(w)$ which is $W^*$ valued and $\nabla u Q$ is trace class on $H$.
Set
$(\trace)_e K=\limn \sum_1^n (e_i, Ke_i)_H$.

\noindent
Then (a) $(\trace)_e (\nabla u Q)$ exists and
\begin{equation}
\label{5.1}
\W\Bigl(\calL_Q w, \tilde u(w)\Bigr)_{W^*} = 
\delta (Qu) + (\trace)_e
(\nabla u q)
\end{equation}
(b)~ If we also assume that $u=\nabla F$ and $Q$ is skew-symmetric then
\begin{equation}
\label{5.2}
_W(\calL_Q w, \widetilde{\widetilde{\nabla F}})_{W^*} = \calL_Q F
\end{equation}
(and then $\calL_Q w$ acts as a vector field on $F$ with $\calL_Q F$ being the
directional derivative along the tangent process).
\end{lemma}

\proof
Setting $(u,e_i) = v_i$
and $u_n = \sum_1^n v_i e_i$, where
$\{e_i\}$ is a CONB induced by $W^*$.
\begin{align*}
\W\Bigl(\calL_\theta w, \tilde{\tilde u}(w) \Bigr)_{W^*} & = 
\sum_i^n \delta (Q e_i) v_i\\
& = \sum_i^n \delta (Q v_i e_i) + \sum_i^n (\nabla v_i, Qe_i) \\
& = 
\calL_Q u_n + \sum_i^n \sum_j^n \nabla_{e_j} v_i (e_j, Qe_i)
\end{align*}
The left hand side is a continuous functional on $W^*$ and the last equation
converges on both sides to \req{5.1}
which proves (a).  (b)~ follows since
$\nabla_{e_j} u_i = \nabla_{e_i, e_j}^2 F$ is symmetric.
\qed

Note that by Lemma 2.1, if
$\calL_Q w$ exist in $\DD_{2,1} (H)$, and if $Q+Q^T=0$ then:
$$
\delta (\calL_Q w) = 0. 
$$

The tangent processes that were considered in \cite{CC}--\cite{D95}  
were of the form of the
right hand side of \req{1.1} in the 
introduction with $\{\sigma_{ij}\}$ skew symmetric and nonanticipative.  
The relation to the
$\calL_Q w$ formulation will now be pointed out.
Consider the case of the $d$-dimensional Brownian motion,
$h=\int_0^\bfcdot h'_s ds$,
$h'\in L^2 ([0,1], \reals^d) $, then we have 
\begin{proposition}\mbox{}\\
(A)~Let $ \boldq$ denote the matrix $\{ q_{ij} (\theta, u), i \le i, j \le d,
\theta, u \in [0,1]\}$ and set
$$
(Q h)_i = \sum_{j=0}^n \int_0^\bfcdot \int_0^1 q_{ij} (\theta, u) h'_j (u) du
d \theta, \quad i=1, \cdots, d
\,.
$$
Assume that the $q_{i,j} (\theta, w)$ are $\calF_\theta$ adapted for all
$u \in [0,1]$ and $E|Qh|_H^2 \le K\cdot| h|_H^2$.

\noindent
Then $\calL_Q w$ exists and as a $H$-valued r.v.
\begin{equation}
\label{nnewstar}
\calL_Q w = \int_0^\bfcdot \left(\int_0^1 \boldq(\theta, u) du \right) 
dw_\theta
\end{equation}
(B)~If $\boldb = \{ b_{i,j} (s), 1\le i, j \le d, s \in [0,1]\}$
where $b_{i,j} (s) $ are $\calF_s$ adapted and
$$
E \sum_{i,j} b_{i,j}^2 (s) < \infty
\,.
$$
Setting
$$
(Bh)_i = \int_0^\bfcdot b_{ij} (s) h'_j(s) ds
$$ 
then
\begin{equation}
\label{nnewstarstar}
\calL_B w = \int_0^\bfcdot \boldb (s) dw_s\,.
\end{equation}
\end{proposition}

\proof Let $(e_i,i\geq 1)$ be an orthonormal basis  of $H$.  Then 
\[
Qe_i=\int_0^\cdot
\left(\int_0^1{\bf{q}}(\tau,u){\dot{e}}_i(u)du\right)d\tau\,,
\]
We have
\begin{align}
\label{sam}
\sum_{i=1}^n \delta(Qe_i)\langle e_i\,, \alpha\rangle
&= \sum_{i=1}^n \int_0^1 
\left(\int_0^1{\bf{q}}(\tau,u){\dot{e}}_i(u)du\right)dw_\tau 
\int_0^1 \dot e_i(s)\dot \alpha(s)ds\nonumber\\
&=  \int_0^1 
\int_0^1{\bf{q}}(\tau,u)\left(\sum_{i=1}^n \int_0^1 \dot e_i(s)\dot \alpha(s)ds
 {\dot{e}}_i(u)\right)  du dw_\tau  \nonumber\\
&= \int_0^1 \int_0^1 {\bf q}(\tau, u)\beta_n(u)du dw_\tau\,,
\end{align}
where \[
\beta_n(u)=\sum_{i=1}^n \int_0^1 \dot e_i(s)\dot \alpha(s)ds\ 
 {\dot{e}}_i(u)\,.
\]
It is clear that $\beta_n$ converges to $\dot \al$ in $L^2[0,1]$.
Since $A$ is of Hilbert-Schmidt, we see that 
\[
\int_0^1 \boldq (\tau, u)\beta_n(u)du
\]
converges to $\int_0^1 \bolda(\tau, u)\dot \alpha(u)du$
in $L^2[0,1]$.  
Consequently, (\ref{sam}) converges to 
\[
\int_0^1 \int_0^1  \boldq(\tau, u)\dot \alpha(u) du dw_\tau\,
\]
in $L^2$ hence in probability and \req{sam} follows.  

In order to prove \req{nnewstarstar} we have to show that
\req{nnewstar} holds for the case where
$$
\boldq (\theta, u) = \bf b (\theta) \delta (\theta -u)
$$
where $\delta$ is the Dirac delta function.
Setting
\begin{align*}
\delta_\eps (u) & = \frac{1}{\eps}u \in [0,\eps] \\
& = 0 \qquad \text{otherwise}
\end{align*}
then, \req{nnewstarstar} follows since by \req{nnewstar} 
$$
\calL_{Q_n} w = \int_0^\bfcdot \frac{1}{\eps} \left(\int_{\theta-\eps}^\theta
b(u) du \right) d w_\theta
$$
and
$$
E \int_0^1 \left(b - \frac{1}{\eps} \int_{\theta-\eps}^\theta b(u) du\right)^2
d\theta \underset{\eps\to0}{\longrightarrow} 0
\,.
$$
\qed

\section{Groups of rotations}
\setcounter{equation}{0}

{\bf A.  Theorem 6.1 (\cite{T}, \cite{G}):}
{\em
Let $T$ be an invertible measure preserving transformation on the Wiener
space,
then there exists a family $T_\theta, \theta \in [0,1]$ of 
measure preserving
transformations such that $T_0 w=w$, $T_1=T$ and for every
$\theta\in[0,1]$,
$E \Bigl|T_\eta w - T_\theta w\Bigr|\W \to0$ as $\eta\to\theta$.
}

The proof was first shown to us by Tsirelson \cite{T}, the proof given here is a
shorter proof due to Glasner \cite{G}.

\proof
Since the Wiener measure on $C[0,1]$ and the Lebesgue measure on $[0,1]$ are
isomorphic, it suffices to prove the result for the Lebesgue measure.  
Set $T_a(X) = a T(X/a)$ in $(0,a)$ and $T_a(X)=X$ in 
$(a,1)$ which is measure preserving,
$T_1=T$ and $T_0$ is the identity.
Now,
\begin{align*}
E\Bigl|T_a(X) - T_{a+\eps}(X)\Bigr| & \le E \left\{\left|
aT\left(\frac{X}{a}\right) - (a+\eps) 
T \left(\frac{X}{a+\eps}\right)\right| \cdot \won_{X \le a}\right\} + a \\
& \le E \left\{\left| a T\left(\frac{X}{a}\right) - 
aT\left(\frac{X}{a+\eps}\right)
\right| \cdot \won_{X \le a}\right\} + 2\eps\,. 
\end{align*}
Applying Lusin's theorem to approximate $T(X)$ by a continuous 
$\tau_\theta(X)$, yields
$$
\left| E\Bigl( Ta(X) - T_{a+\eps} (X)\Bigr) \right| \le 
a\left| E \left( \tau_\theta \left(\frac{X}{a}\right) - \tau_\theta 
\left(\frac{X}{a=\eps}\right)
\right) \right| + 2\theta + 2 \eps
$$
and continuity in $\eps$ follows by dominated convergence since $\theta$ is
arbitrary. 
\qed

\vspace{-1cm}
\paragraph{B.  Flows}\mbox{}\\
We want to show that for $A_t+A_t^T = 0$ and additional conditions
the equation
\begin{equation}
\label{6.1}
\frac{dT_tw}{dt} = \Bigl( \calL_{A_t(w)} w\Bigr) \circ T_tw, \quad T_0w=w
\end{equation}
defines a flow of rotations.
The case where $W$ is the $d$-dimension Wiener space and $A$ is adapted:
$$
\Bigl( \calL_A w\Bigr)_i = \int_0^\bfcdot \sum_j^d a_{ij} (s,w) dw_j(s),
\qquad i=1,2, \cdots, d
$$
was considered by Cipriano and Cruzeiro \cite{CC}.  
The general 
result presented in the next theorem is from \cite{HUZ} and is followed
by a more detailed proof.
%%%%%%%%%%%%%%%%%%%%%%%%%%%%%%%%%%%%%%%%%%%%%%%%%%%%%%%%%
\setcounter{theorem}{1}
\begin{theorem}[\cite{HUZ}]
Assume that for all $t\ge 0$, 
$A=A_t (w) :$ is a skew symmetric  strongly measurable
mapping and  for  any $h\in H$, $A_th\in \DD_{p,1} (H) $,
for some $p>1$ and a.a.\ $t\in [0,T]$
where $T>0$ is fixed. 
Further assume that:
\begin{enumerate}
\item The series 
$$
B_t=\sum_{i=1}^\infty \delta(A_te_i)e_i
$$
converges 
in $L^p(d\mu\times dt,W)$ (as a $W$-valued random
variable), where $(e_i,i\geq 1)$ is a fixed orthonormal basis of $H$.
\item Let $p_n$ denote the orthogonal
  projection onto the span of $\{e_1,\ldots,e_n\}$ and $V_n$ denote the
  sigma-algebra generated by $\{\delta e_1,\ldots,\delta e_n\}$.  
Assume that the sequence of vector fields $(B^n,n\geq 1)$ defined by
\begin{equation}
\label{6.2}
B_t^n = \sum_{i=1}^n\delta E[p_n A_t e_i|V_n]e_i
\end{equation}
or, as will be shown later to be the same as
\begin{equation}
\label{6.2a}
B_t^n =\sum_{i=1}^n E[\delta( p_n A_t e_i)|V_n]e_i \tag{6.2a}
\end{equation}
converges to $B$ in $L^p(d\mu\times dt,W)$.
\item 
Assume that
for a given $\eps>0$, we have 
\begin{equation}
\label{6.3}
\int_0^TE\left\{\exp\eps\||\nabla B_t\||\right\}dt=\Gamma_{0,T} <\infty\,,
\end{equation}
the norm above is defined  as 
$$
\||\nabla B_t\||=\sup\{|\,\nabla_h(_W\lip B_t, \alpha\rip_{W^*})\,|\,:h\in
B_1,\,\alpha\in W^\star_1\}
$$
where $B_1 = \{ h \in H: |h|_W = 1\}$ and
$W_1^\star $ is the unit ball of $W^\star$.
Further assume that \req{6.3} also holds for $B_t$ replaced by $B_t^n$ (this
holds, e.g., when $\{e_i\}$ is a Schauder basis of $W$).
\end{enumerate}
Then the equation 
\begin{equation}
\label{flow-1}
\phi_{s,t}(w)=w+\int_s^tB_r(\phi_{s,r})dr,\,\,s<t\,,
\end{equation}
defines a flow of measure preserving diffeoemorphisms of $W$  whose
almost sure inverse is denoted by 
$(\psi_{s,t},\,0\leq s\leq t\leq T)$ and  satisfies 
$$
\mu\{w:\,\phi_{s,t}\circ\psi_{s,t}(w)=\psi_{s,t}\circ\phi_{s,t}(w)=w\}=1\,.
$$
Moreover the inverse flow is the unique solution of the 
equation 
\begin{equation}
\label{flow-2}
\psi_{s,t}(w)=w-\int_s^tB_r(\psi_{r,t})dr
\end{equation}
and $\phi_{s,t}$ (hence $\psi_{s,t}$) leaves the Wiener measure
invariant, i.e. $\phi_{s,t}^\star\mu=\mu$  for any $s<t\in [0,T]$.
\end{theorem}

\proof
We start with showing the equality
\req{6.2} and \req{6.2a}.
Set $\alpha^m = p_m Ap_m$ and $a_{ij}^m = (e_j, \alpha^m e_i)$, then for 
$i, j \le m $
\begin{align*}
E\Bigl\{\delta (p_m A e_i) | V_m\Bigr\} = 
E\Bigl\{\delta (\alpha^m e_i) | V_m\Bigr\} 
& = \sum_{j=1}^m E \Bigl\{ \delta (a_{ij}^m e_j) | V_m\Bigr\} \\
& = \sum_{j=1}^m \delta e_j E\{ a_{i,j} | V_m\} - \sum_j E \{ \nabla _{e_j}
a_{ij} | V_m\} \\
& = 
\sum_{j=1}^m \delta e_j E \Bigl\{ a_{ij} | V_m\Bigr\} 
- \sum_j \nabla_{e_j} E \{a_{ij} |
V_m\}\\
& = \sum_{j=1}^m \delta \left(E\Bigl\{a_{ij} | V_m\Bigr\} \cdot e_j\right) \\
& = \delta E(\alpha^m e_i | V_m)
\end{align*}
and \req{6.2a} follows. 
%Consequently, \req{6.3} implies that
%$$
%\sup_{m} \int_0^T E \exp \eps |\| \nabla B_t^m |\| dt < \infty\,.
%$$
Also, since $A_t$ is skew symmetric, so are the matrices
$a_{ij}^m$ and $E(a_{ij}^m | V_m)$ hence
$\delta B_t^m = 0$.  Consequently, (cf., e.g., \cite[Theorem~5.3.1]{UZ99}
the claimed results of Theorem~6.1 hold for $B_t$ replaced by $B_t^m$.
Therefore, denoting by
$\phi_{s,t}^n, s \le t\in [0,T]$) the flow associated to the cylindrical
vector field $B^n$ and by $\psi_{s,t}^n, s \le t \in [0,T]$ its inverse, then
in particular we have
\begin{align*}
\frac{d \phi_{s,t}^{n\star} \mu}{d\mu} & = \exp \int_s^t
\delta B_r^n (\psi_{r,t}^n) dr 
 = 1
\end{align*}
and
\begin{align*}
\frac{d \psi_{s,t}^{n\star} \mu}{d\mu} & = \exp - \int_s^t
\delta B_r^n (\phi_{s,r}^n) dr 
 = 1\,.
\end{align*}

Let $e_1,e_2, \cdots $ be a fixed CONB of $H$ induced by elements of $W^*$.
Let $M$ denote the following class of cylindrical operator $Q$ on $H$.  Let
$q_{ij} = (e_j, Q e_i)$ then, for some $m$\\
(a)~$q_{i,j} = 0$ for $i>m$ or $j>m$\\
(b)~$q_{i,j} = -q_{ji}$\\
(c)~$q_{ij} = f_{ij}(\delta e_1, \cdots, \delta e_m)$ and $f_{ij}$ possesses
bounded first derivatives.

Set $B_r(w) = \sum_1^m \delta (Q_r (w) e_i) e_i, Q_r \in M$.
The following version of \cite[Theorem~5.2.1]{UZ99}
is needed to complete the  proof of the theorem.
\begin{proposition}
Let $Q_r^a, Q_r^b \in M$, $r\in [0,T]$ and assume
that for some $\eps > 0$
\begin{equation}
\label{6.6}
E \int_0^T \Bigl(\exp \eps |\| \nabla B_r^a |\| + 
\exp \eps |\| \nabla B_r^b|\| \Bigr)
dr \le \Gamma_{0,T} < \infty\,.
\end{equation}
Let $\varphi_{s,t}^a, \varphi_{s,t}^b$ denote the flows induced by $Q^a$ and
$Q^b$.  Then for $s<t$, $(t-s)$ sufficiently small
$$
E \sup_{u\in[s,t]} \Bigl| \varphi_{s,u}^a - \varphi_{s,u}^b \Bigr|\W
\le E \left( \int_s^t \Bigl| B_r^a - B_r^b \Bigr|\W^p dr\right)^{\frac{1}{p}}
\Gamma_{s,t}^{\frac{1}{\gamma}} \left(\frac{1}{t-s}\right)^{\frac{1}{q}}
$$
where $\Gamma_{s,t}$ is defined as $\Gamma_{0,T}$ (equation \req{6.3}) with
$0,T$ replaced by $s,t,$ and $(t-s) q \le \eps$.
\end{proposition}
\paragraph{Proof of proposition:}
Set $D_r = B_r^a - B_r^b$ and let $\varphi_{s,t}^\lambda$,
$\la \in [0,1]$ be the solution to
$$
\varphi_{s,t}^\la (w) = w+ \int_s^t \Bigl(\la B_r^a + (1-\la) B_r^b\Bigr)
\circ \varphi_{s,r}^\la dr\,.
$$
Then $\varphi_{s,t}^\la$ is also a rotation.  Set 
$Z_{s,t}^\la = \frac{d \varphi_{s,t}^\la}{d\la}$, then
$$
\varphi_{s,t}^a - \varphi_{s,t}^b = \int_0^1 Z_{s,t}^\la d \la
$$
and
$$
Z_{s,t}^\la = \int_s^t D_r \circ \varphi_{s,r}^\la d \la + \int_s^t \left[(\nabla
(B_r^b + \la D_r)\right] \circ \varphi_{s,r}^\la Z_{s,r}^\la dr
\,.
$$
By Gronwall's lemma
$$
\Bigl|Z_{s,t}^\la\Bigr|\W \le \int_s^t \Bigl|D_r \circ \varphi_{s,r}^\la 
\Bigr|\W dr \cdot
\exp \int_s^t |\| \nabla B_r^b + \la \nabla D_r |\| \circ \varphi_{s,r}^\la dr
\,.
$$
Therefore, since 
$\varphi^\la$ is measure preserving:
\begin{align*}
E & \sup_{r\in [s,t]}  \Bigl| \varphi_{s,r}^a - \varphi_{s,r}^b\Bigr|\W \\
& \le E \int_0^1 \Bigl| Z_{s,t}^\la \Bigr|\W d\la \\
& \le E \left\{ \int_0^1 \left( \int_s^t \Bigl| D_r (w) \Bigr|\W dr\right)
\exp \int_s^t |\| \nabla B_r^b + \la \nabla D_r |\| dr d\la \right\} \\
& \le E \left\{ \int_s^t \Bigl|D_r\Bigr|\W dr \exp \int_s^t \Bigl( |\| \nabla B_r^b|\| +
|\| \nabla B_r^a|\|\Bigr) dr \right\} \\
& \le \left(E\int_s^t \Bigl|D_r\Bigr|\W^p dr\right)^{\frac{1}{p}}
E\left( \exp\: q\int_s^t \Bigl( |\| \nabla B_r^a |\| + |\| \nabla B_r^a|\|\Bigr)
dr \right)^{\frac{1}{q}} \\
&\le \left(E \int_s^t \Bigl|D_r\Bigr|\W^P dr \right)^{\frac{1}{p}}
E \left(\frac{1}{t-s} \int_s^t \exp \: q(t-s) \Bigl(|\| \nabla B_r^a |\| +
|\| B_r^b|\|\Bigr) dr \right)^{\frac{1}{q}}
\end{align*}
which proves the proposition.
\qed

Returning to the proof of the theorem, setting
$$
Q_r^a = E(p_m A_r p_m | V_m)
$$
and similarly, with $m$ replaced by $n$, for $Q_r^b$ yields 
\begin{equation}
\label{6.8}
E\sup_{u\in[s,t]} \Bigl| \varphi_{s,u}^m - \varphi_{s,u}^n\Bigr|\W
\le \left(\Gamma_{s,t} \frac{1}{t-s}\right)^{\frac{1}{q}}
E \int_s^t \Bigl|B_r^m - B_r^n\Bigr|\W dr\,.
\end{equation}

This result implies the  convergence of
$(\phi^n_{s,u},\,u\in [s,t])$ in $L^1(\mu,W)$  uniformly with respect
to $u$ for the intervals $[s,t]$ and the limit $(\phi_{s,u},u\in
[s,t])$ is the unique solution of the equation (\ref{flow-2}). Equation
(\ref{6.8}) implies the uniqueness of the since if
$(\phi'_{s,u})$ is another solution, then its
finite dimensional approximations must coincide with
$(\phi^n_{s,u})$. Now, using Lemmas 5.3.1,  5.3.2 and  5.3.5
of \cite{UZ99}, it can be shown that the constructions of $(\phi_{s,u})$
on the different small intervals can be patched together to give the
entire flow. For the inverse flow the same  reasoning  applies also.
\qed

\section*{Appendix}
\renewcommand{\theequation}{A.\arabic{equation}}
    \setcounter{equation}{0}

\hsp
As discussed in section~2, a necessary condition for $u\in \DD_{2,1}(H)$ to be
adapted to a given filtration is that $\nabla u$ be a.s. quasinilpotent.
Now, given a $u\in \DD_{2,1}(H)$ such that $\nabla u$ is quasinilpotent,
the question arises whether this assumes the existence of a filtration for
which $u$ is adapted.  The answer is negative as seen from the following
example:

Let
$$
u(w) = \sum_{i=1}^\infty 2^{-i\alpha} \delta (e_{i+1}) e_i
$$
where $e_i$ is a CONB in $H$.  Hence, for $h\in H$
\begin{align*}
(\nabla u) h &= \sum_{i=1}^\infty 2^{-i\alpha} (e_{iH}, h) e_i \\
(\nabla u)^r h & = \sum_{i=1}^\infty \beta_{i,r} (e_{i+r}, h) e_i
\end{align*}
where
$\beta_{i,r}  = 2^{-i \alpha}$,  
$\beta_{i+1, r-1} = 2^{-\alpha r(2i + r + 1)/2}$.
Following Ringrose \cite{R}:
\begin{enumerate}
\item
$\nabla R$ is quasinilpotent.
\item
For any $h\in H$, it holds for $r$ large enough that
\begin{equation}
\label{Astar}
\left|\Bigl((\nabla u)^r h, e_i\Bigr)_H \left| > \half \right|(\nabla u)^r
h\right|_H\,.
\end{equation}
\end{enumerate}

\renewcommand{\thelemma}{A}
\setcounter{theorem}{0}

\begin{lemma}
Let $E_\bfcdot$
be a strictly monotone resolution of the identity.  Given $\delta>0$ and some
$e$ in $H$, then there
exists an $N$ such that for all $X$ for which $(I-E_{1-\frac{1}{N}}) X=X$,
$|X|_H =1$, we have $|(X,e)_H|\le \delta$.
\end{lemma}

\proof
Assume that the Lemma is not true, then for every $N$ there exists
an $X$ satisfying $(I-E_{1-\frac{1}{N}}) X=X$,
$|X|_H=1$ and $|(X,e)| > \delta$.  Let $M_N$ denote the collection of such
$X's$.  Then, by compactness, for every $N$ there is a converging sequence and
then consider a collection of such sequences as $N\to\infty$
from which we can construct a converging sequence with $N\to\infty$.
Let $X_n$ $N\to\infty$ denote this subsequence then $|X_W|_H = 1$,
$X_N = (I-E_{1-\frac{1}{N}}) X_N$.
Hence the limit must vanish since $E_\bfcdot$ is continuous contradicts our
assumption.
\qed

Assume now that $u$ is adapted to $\calF_\bfcdot$ induced by 
$E_\bfcdot$ then 
$$
\nabla u \cdot E_\la = E_\la \cdot \nabla u \cdot E_\la
\,.
$$
Hence
$$
(\nabla \,u)^T (I-E_\la) = (I-E_\la) (\nabla\, u)^T (I-E_\la)\,.
$$
Consequently for $\la=-\frac{1}{N}$, if $X\in M_N$ so in $(\nabla u)^*X$.
Hence, for $r$ large enough by \req{Astar}
$$
\left|\Bigl( (\nabla\,u)^T\Bigr)^r X, e_i \right| \ge \half
$$
and by the lemma the left hand side is upper bounded by any $\delta < \half$
hence $u$ cannot be adapted to any continuous filtration.

\noindent
M. Zakai\\
Department of Electrical Engineering\\
Technion--Israel Institute of Technology\\
Haifa 32000, Israel\\
{\tt zakai@ee.technion.ac.il}

\end{document}